\documentclass[11pt,a4paper]{amsart}
\usepackage{amsthm,amsmath,amsfonts,calc,graphicx}
\usepackage[sort&compress,numbers]{natbib}
\usepackage[lmargin=40mm,rmargin=40mm,tmargin=35mm,bmargin=35mm]{geometry}

\theoremstyle{plain}
\newtheorem{theorem}{Theorem}[section]
\newtheorem{lemma}[theorem]{Lemma}
\newtheorem{corollary}[theorem]{Corollary}

\newcommand{\comment}[1]{
\medskip\noindent
\framebox{\parbox{\textwidth-2mm}{#1}}
\medskip}

\renewcommand{\comment}[1]{}

\newcommand{\spacing}[1]{
\renewcommand{\baselinestretch}{#1}
\setlength{\footnotesep}{\baselinestretch\footnotesep}}
\spacing{1.15}

\newcommand{\corref}[1]{Corollary~\ref{cor:#1}}

\newcommand{\twothmref}[2]{Theorems~\ref{thm:#1} and \ref{thm:#2}}
\newcommand{\thmref}[1]{Theorem~\ref{thm:#1}}
\newcommand{\lemref}[1]{Lemma~\ref{lem:#1}}
\newcommand{\twolemref}[2]{Lemmata~\ref{lem:#1} and \ref{lem:#2}}
\newcommand{\figref}[1]{Figure~\ref{fig:#1}}
\newcommand{\secref}[1]{Section~\ref{sec:#1}}
\newcommand{\SECREF}[1]{\S\ref{sec:#1}}
\newcommand{\seclabel}[1]{\label{sec:#1}}
\newcommand{\figlabel}[1]{\label{fig:#1}}

\newcommand{\thmlabel}[1]{\label{thm:#1}}

\newcommand{\lemlabel}[1]{\label{lem:#1}}
\newcommand{\corlabel}[1]{\label{cor:#1}}
\newcommand{\mySection}[2]{\section{#1}\seclabel{#2}}
\newcommand{\twosecref}[2]{Sections~\ref{sec:#1} and \ref{sec:#2}}

\newcommand{\TB}[1]{\ensuremath{\protect\Theta(#1)}}
\newcommand{\Oh}[1]{\ensuremath{\protect\mathcal{O}(#1)}}
\newcommand{\arc}[1]{\ensuremath{\protect\overrightarrow{#1}}}
\newcommand{\etal}{~et~al.~}
\newcommand{\half}{\ensuremath{\protect\tfrac{1}{2}}}
\newcommand{\third}{\ensuremath{\protect\tfrac{1}{3}}}
\newcommand{\floor}[1]{\ensuremath{\protect\lfloor#1\rfloor}}
\newcommand{\ceil}[1]{\ensuremath{\protect\lceil#1\rceil}}

\newcommand{\aaa}{\textup{\hspace*{0.5em}(a)\hspace*{0.25em}}}
\newcommand{\bbb}{\textup{\hspace*{0.5em}(b)\hspace*{0.25em}}}
\newcommand{\ccc}{\textup{\hspace*{0.5em}(c)\hspace*{0.25em}}}

\newcommand{\Z}{\ensuremath{\mathbb{Z}}}
\newcommand{\R}{\ensuremath{\mathbb{R}}}

\newcommand{\x}{\ensuremath{\protect\textup{\textsf{x}}}}
\newcommand{\y}{\ensuremath{\protect\textup{\textsf{y}}}}
\newcommand{\z}{\ensuremath{\protect\textup{\textsf{z}}}}

\newcommand{\cn}[1]{\ensuremath{\chi(#1)}}
\newcommand{\sscn}[1]{\ensuremath{\chi_{\textup{sst}}(#1)}}
\newcommand{\tn}[2][]{\ensuremath{\textup{\textsf{tn}}_{#1}(#2)}}
\newcommand{\utn}[2][]{\ensuremath{\textup{\textsf{utn}}_{#1}(#2)}}
\newcommand{\qn}[2][]{\ensuremath{\textup{\textsf{qn}}_{#1}(#2)}}
\newcommand{\uqn}[2][]{\ensuremath{\textup{\textsf{uqn}}_{#1}(#2)}}

\newcommand{\Figure}[4][htb]{
\begin{figure}[#1]
	\vspace*{1ex}
	\begin{center}#3\end{center}
	\vspace*{-1ex}
	\caption{\figlabel{#2}#4}
\end{figure}
}

\begin{document}

\title[Upward Three-Dimensional Grid Drawings of Graphs]{\textbf{Upward Three-Dimensional\\ Grid Drawings of Graphs}}

\author{Vida Dujmovi{\'c}}
\address{School of Computer Science, Carleton University, Ottawa, Canada}
\email{vida@scs.carleton.ca} 

\author{David R. Wood}
\address{Departament de Matem{\`a}tica Aplicada II, Universitat Polit{\`e}cnica de Catalunya, Barcelona, Spain}
\email{david.wood@upc.edu}


\thanks{Research of Vida Dujmovi\`c is supported by NSERC. Research of David Wood is supported by the Government of Spain grant MEC SB2003-0270 and by the projects MCYT-FEDER BFM2003-00368 and Gen.\ Cat 2001SGR00224.}


\begin{abstract}  
A \emph{three-dimensional grid drawing} of a graph is a placement of the vertices at distinct points with integer coordinates, such that the straight line segments representing the edges do not cross. Our aim is to produce three-dimensional grid drawings with small bounding box volume. Our first main result is that every $n$-vertex graph with bounded degeneracy has a three-dimensional grid drawing with \Oh{n^{3/2}} volume. This is the largest known class of graphs that have such drawings. A three-dimensional grid drawing of a directed acyclic graph (\emph{dag}) is \emph{upward} if every arc points up in the \z-direction. We prove that every dag has an upward three-dimensional grid drawing with \Oh{n^3} volume, which is tight for the complete dag. The previous best upper bound was \Oh{n^4}. Our main result concerning upward drawings is that every $c$-colourable dag ($c$ constant) has an upward three-dimensional grid drawing with \Oh{n^2} volume. This result matches the bound in the undirected case, and improves the best known bound from \Oh{n^3} for many classes of dags, including planar, series parallel, and outerplanar.  Improved bounds are also obtained for tree dags. We prove a strong relationship between upward three-dimensional grid drawings, upward track layouts, and upward queue layouts. Finally, we study upward three-dimensional grid drawings with bends in the edges. 
\end{abstract}

\maketitle

%%%%%%%%%%%%%%%%%%%%%%%%%%%%%%%%%%%%%%%%%%%%%%%%%%%%%%%%%%%%%%%%%%%%%%%%%%%%%%
\mySection{Introduction}{Introduction}
%%%%%%%%%%%%%%%%%%%%%%%%%%%%%%%%%%%%%%%%%%%%%%%%%%%%%%%%%%%%%%%%%%%%%%%%%%%%%%

Graph drawing is the study of aesthetically pleasing geometric representations of graphs. Graph drawing in the plane is well-studied; see \cite{DETT99, KaufmannWagner01}. Motivated by experimental evidence suggesting that displaying a graph in three dimensions is better than in two \cite{WF94, WF96}, and applications  including information visualisation \cite{WF94}, VLSI circuit design \cite{LR86}, and software engineering \cite{WHF93}, there is a growing body of research in three-dimensional graph drawing. 

A \emph{three-dimensional straight line grid drawing} of a graph, henceforth called a \emph{3D drawing}, is a placement of the vertices at distinct points in $\mathbb{Z}^3$ (called \emph{gridpoints}), such that the straight line segments representing the edges are pairwise non-crossing. That is, distinct edges only intersect at common endpoints, and each edge only intersects a vertex that is an endpoint of that edge. The coordinates of a vertex $v$ are denoted by $(\x(v),\y(v),\z(v))$. It is well known that every graph has a 3D drawing.  We are therefore interested in optimising certain measures of the aesthetic quality of such drawings. 

The \emph{bounding box} of a 3D drawing  is the minimum axis-aligned box that contains the drawing. If the bounding box has side lengths $X-1$, $Y-1$ and $Z-1$, then  we speak of an $X\times Y\times Z$ \emph{drawing} with \emph{width} $X$, \emph{depth} $Y$, \emph{height} $Z$,  and \emph{volume} $X\cdot Y\cdot Z$. That is, the volume of a 3D drawing is the number of gridpoints in the bounding box. This definition is formulated so that 2D drawings have positive volume. We are interested in 3D drawings with small volume, which are widely studied \citep{BCMW-JGAA04, CS-IPL97, CELR-Algo96, DMW-SJC05, FLW-JGAA03, PTT99, Poranen-00, GLM-CGTA05, DujWoo-DMTCS05, DujWoo-SubQuad-AMS, Giacomo-GD03, dGLW-CCCG02, DM-GD03, Hasunuma-GD03}. 

3D drawings have been generalised in a number of ways. Multi-dimensional grid drawings have been studied \citep{PorWoo-GD04, Wood-Hypercube}, as have 3D \emph{polyline} grid drawings, where edges are allowed to bend at gridpoints \citep{DELPW-GD05, DujWoo-DMTCS05, Wismath-TR04, MorinWood-JGAA04}. The focus of this paper is upward 3D drawings of directed graphs, which have previously been studied by \citet{Poranen-00} and \citet{DLMW-GD05}. A 3D drawing of a directed graph $G$ is \emph{upward} if $\z(v)<\z(w)$ for every arc \arc{vw} of $G$. Obviously an upward 3D drawing can only exist if $G$ is acyclic (a \emph{dag}). Upward two-dimensional drawings have been widely studied; see \citep{BDMT98, GT-SJC01, HL-SJC96, KR-Order93, Thom-Order89, GT-ORDER95} for example.

As described in Table~\ref{tab:VolumeUpperBounds}, our main results are improved upper bounds on the volume of upward 3D drawings of dags. These results are presented in \twosecref{Arbitrary}{Coloured}, and in \secref{Trees} in the case of trees.  In addition, we prove that (undirected) graphs with bounded degeneracy have 3D drawings with \Oh{n^{3/2}} volume. This is the largest known class of graphs that have such drawings.

\begin{table}[htb]
\begin{center}
\caption{\label{tab:VolumeUpperBounds}
Upper bounds on the volume of 3D drawings of graphs and upward 3D drawings of dags with $n$ vertices, $m$ edges, chromatic number $\chi$, and degeneracy $d$.}
\begin{tabular}{|l|lr|lr|lr|}
\hline
&& 
& \multicolumn{4}{c|}{upward dags}\\\cline{4-7}
\raisebox{1.5ex}[0cm][0cm]{graph family}
& \multicolumn{2}{c|}{\raisebox{1.5ex}[0cm][0cm]{undirected}}
& \multicolumn{2}{c|}{previous best}
& \multicolumn{2}{c|}{this paper}
\\\hline
%%%%%%%%%%%%%%%%%%%%%%%%%%%%%%%%%%%%%%%%%%%%%%%%%%%%%%%%%%%
arbitrary
	& \TB{n^3}	& \citep{CELR-Algo96}	
	& \Oh{n^4}	& \citep{DLMW-GD05}		
	& \TB{n^3}	& \SECREF{Arbitrary}\\
%%%%%%%%%%%%%%%%%%%%%%%%%%%%%%%%%%%%%%%%%%%%%%%%%%%%%%%%%%%
arbitrary					
	& \Oh{m^{4/3}n}	& \citep{DujWoo-SubQuad-AMS}
	& & & & \\
%%%%%%%%%%%%%%%%%%%%%%%%%%%%%%%%%%%%%%%%%%%%%%%%%%%%%%%%%%%
arbitrary 							
	& \Oh{dmn}		& \SECREF{Strong}
	& & & & \\
%%%%%%%%%%%%%%%%%%%%%%%%%%%%%%%%%%%%%%%%%%%%%%%%%%%%%%%%%%%
arbitrary 							
	& \Oh{\chi^2n^2}		& \citep{PTT99}
	& 			&
	& \Oh{\chi^4n^2}		& \SECREF{Coloured}\\
%%%%%%%%%%%%%%%%%%%%%%%%%%%%%%%%%%%%%%%%%%%%%%%%%%%%%%%%%%%
constant $\chi$			
	& \TB{n^2}		& \citep{PTT99}
	& 				& 
	& \TB{n^2}		& \SECREF{Coloured}\\
%%%%%%%%%%%%%%%%%%%%%%%%%%%%%%%%%%%%%%%%%%%%%%%%%%%%%%%%%%%
constant $\chi$			
	& \Oh{m^{2/3}n}	& \citep{DujWoo-SubQuad-AMS}
	& 				& 
	& \Oh{n^2}		& \SECREF{Coloured}\\
%%%%%%%%%%%%%%%%%%%%%%%%%%%%%%%%%%%%%%%%%%%%%%%%%%%%%%%%%%%
%no $H$-minor		
minor-closed 
	& \Oh{n^{3/2}}	& \citep{DujWoo-SubQuad-AMS}
	& 				& 
	& \Oh{n^2}		& \SECREF{Coloured}\\
%%%%%%%%%%%%%%%%%%%%%%%%%%%%%%%%%%%%%%%%%%%%%%%%%%%%%%%%%%%
constant $d$
	& \Oh{n^{3/2}}	& \SECREF{Strong}
	& 				& 
	& \Oh{n^2}		& \SECREF{Coloured}\\
%%%%%%%%%%%%%%%%%%%%%%%%%%%%%%%%%%%%%%%%%%%%%%%%%%%%%%%%%%%
planar		 			
	& \Oh{n^{3/2}}	& \citep{DujWoo-SubQuad-AMS}
	& \Oh{n^3}		& \citep{DLMW-GD05}	
	& \Oh{n^2}		& \SECREF{Coloured}\\
%%%%%%%%%%%%%%%%%%%%%%%%%%%%%%%%%%%%%%%%%%%%%%%%%%%%%%%%%%%
constant treewidth		 			
	& \Oh{n}		& \citep{DMW-SJC05}
	& \Oh{n^3}		& \citep{DLMW-GD05}	
	& \Oh{n^2}		& \SECREF{Coloured}\\
%%%%%%%%%%%%%%%%%%%%%%%%%%%%%%%%%%%%%%%%%%%%%%%%%%%%%%%%%%%
series parallel
	& \Oh{n}		& \citep{DMW-SJC05}
	& \Oh{n^3}		& \citep{Poranen-00}
	& \Oh{n^2}		& \SECREF{Coloured}\\
%%%%%%%%%%%%%%%%%%%%%%%%%%%%%%%%%%%%%%%%%%%%%%%%%%%%%%%%%%%
outerplanar		 			
	& $2\times2\times n$	& \citep{FLW-JGAA03}
	& \Oh{n^3}		& \citep{DLMW-GD05}	
	& \Oh{n^2}		& \SECREF{Coloured}\\
%%%%%%%%%%%%%%%%%%%%%%%%%%%%%%%%%%%%%%%%%%%%%%%%%%%%%%%%%%%
trees
	& $2\times2\times n$				& \citep{FLW-JGAA03}
	& $7\times 7\times 7n$				& \citep{DLMW-GD05}
	& $4\times4\times\frac{7}{5}n$		& \SECREF{Trees}\\
%%%%%%%%%%%%%%%%%%%%%%%%%%%%%%%%%%%%%%%%%%%%%%%%%%%%%%%%%%%
caterpillars
	& $2\times2\times n$	& \citep{FLW-JGAA03}
	& 						&
	& $2\times2\times n$	& \SECREF{Trees}\\
\hline
\end{tabular}
\end{center}
\end{table}

Other results in this paper include the following. In \secref{Track} we study upward track layouts, and show how they can be used to produce upward 3D drawings with small volume. These results are used in \secref{Trees} to produce upward 3D drawings of trees. In \secref{Queue} we explore the relationship between upward track layouts and upward queue layouts, which is a structure introduced by Heath\etal\citep{HP-SJC99, HPT-SJC99} in the 1990's. In \secref{Example} we describe an outerplanar graph that highlights the key differences between 3D drawings and upward 3D drawings. Finally in \secref{Subdivisions} we study upward layouts of graph subdivisions, and conclude with some bounds on the volume of upward 3D polyline drawings.

%%%%%%%%%%%%%%%%%%%%%%%%%%%%%%%%%%%%%%%%%%%%%%%%%%%%%%%%%%%%%%%%%%%%%%%%%%%%%%
\section{Preliminaries}
%%%%%%%%%%%%%%%%%%%%%%%%%%%%%%%%%%%%%%%%%%%%%%%%%%%%%%%%%%%%%%%%%%%%%%%%%%%%%%

The following notation is used throughout the paper. We consider finite simple graphs $G$ with vertex set $V(G)$. If $G$ is undirected then its edge set is denoted by $E(G)$. If $G$ is directed then its arc set is denoted by $A(G)$. A \emph{vertex ordering} of $G$ is a bijection $\sigma:V(G)\rightarrow\{1,2,\dots,n\}$, sometimes written as $\sigma=(v_1,v_2,\dots,v_n)$ where $\sigma(v_i)=i$. A vertex ordering $\sigma$ of a directed graph $G$ is \emph{topological} if $\sigma(v)<\sigma(w)$ for every arc $\arc{vw}\in A(G)$. It is well known that a directed graph is acyclic if and only if it has a topological vertex ordering. 

A (\emph{vertex}) $c$-\emph{colouring} of a graph $G$ is a partition $\{V_i:i\in I\}$ of $V(G)$, such that $|I|=c$, and for every edge $vw\in E(G)$, if $v\in V_i$ and $w\in V_j$ then $i\ne j$. Each $i\in I$ is a \emph{colour}, each set $V_i$ is a \emph{colour class}, and if $v\in V_i$ then $v$ is \emph{coloured} $i$. If $G$ has a vertex $c$-colouring then $G$ is $c$-\emph{colourable}. The \emph{chromatic number} of $G$, denoted by \cn{G}, is the minimum integer $c$ such that $G$ is $c$-colourable. 

A graph $G$ is \emph{$d$-degenerate} if every subgraph of $G$ has a vertex of degree at most $d$. The \emph{degeneracy} of $G$ is the minimum integer $d$ such that $G$ is $d$-degenerate. A $d$-degenerate graph is $(d+1)$-colourable by a greedy algorithm. For example, every forest is $1$-degenerate, every outerplanar graph is $2$-degenerate, and every planar graph is $5$-degenerate. 
%In general, degeneracy is at most $\Delta$ and less than $\sqrt{2m}$. 

%%%%%%%%%%%%%%%%%%%%%%%%%%%%%%%%%%%%%%%%%%%%%%%%%%%%%%%%%%%%%%%%%%%%%%%%%%%%%%
\section{Arbitrary Graphs}
\seclabel{Arbitrary}
%%%%%%%%%%%%%%%%%%%%%%%%%%%%%%%%%%%%%%%%%%%%%%%%%%%%%%%%%%%%%%%%%%%%%%%%%%%%%%

\citet{CELR-Algo96} proved that every graph has a 3D drawing with \Oh{n^3} volume. The proof generalises for upward 3D drawings as follows.

\begin{theorem} 
Every dag $G$ on $n$ vertices has an upward $2n\times 2n\times n$ drawing with $4n^3$ volume. Moreover, the bounding box of every upward 3D drawing 
of the complete dag on $n$ vertices is at least $\frac{n}{4}\times\frac{n}{4}\times n$, and thus has 
$\Omega(n^3)$ volume.
\end{theorem}

\begin{proof} 
Let $(v_1,v_2,\dots,v_n)$ be a topological vertex ordering of $G$. By Bertrand's Postulate there is a prime number $p$ such that $n<p\leq2n$. Position each vertex $v_i$ at $(i^3\bmod p,i^2\bmod p,i)$. \citet{CELR-Algo96} proved that no two arcs cross (with the \x- and \z-coordinates switched). Clearly every arc is upward. The bounding box is at most $2n\times 2n\times n$. \citet{CELR-Algo96} observed that the bounding box in every 3D drawing of $K_n$ is at least $\frac{n}{4}\times\frac{n}{4}\times\frac{n}{4}$ (since at most four vertices can lie in a single gridplane). The same lower bound holds for upward 3D drawings of the $n$-vertex complete dag. In addition, the height is at least $n$, since the complete dag contains a Hamiltonian directed path; see \lemref{LongPathLowerBound} below. 
\end{proof}

%%%%%%%%%%%%%%%%%%%%%%%%%%%%%%%%%%%%%%%%%%%%%%%%%%%%%%%%%%%%%%%%%%%%%%%%%%%%%%
\section{Coloured Graphs}
\seclabel{Coloured}
%%%%%%%%%%%%%%%%%%%%%%%%%%%%%%%%%%%%%%%%%%%%%%%%%%%%%%%%%%%%%%%%%%%%%%%%%%%%%%

\citet{PTT99} proved that every $c$-colourable graph has a $\Oh{c}\times\Oh{n}\times\Oh{cn}$ drawing with \Oh{c^2n^2} volume. The proof implicitly relied on the following constructions.

\begin{lemma}[\citep{PTT99}]
\lemlabel{Pach}
Let $\{V_i:0\leq i\leq c-1\}$ be a $c$-colouring of a graph $G$. 
Let $p\geq2c-1$ be a prime number. Place each vertex in $V_i$ at a distinct gridpoint $(i,t,it)$, where $t\equiv i^2\bmod{p}$. Then a \textup{(}crossing-free\textup{)} 3D drawing of $G$ is obtained.
\end{lemma}

\begin{lemma}[\citep{PTT99}]
\lemlabel{PachPach}
Let $\{V_i:0\leq i\leq c-1\}$ be a $c$-colouring of an $n$-vertex graph $G$. 
Then $G$ has a $\Oh{c}\times\Oh{n}\times\Oh{cn}$ drawing, such that $\x(v)<\x(w)$ for all vertices $v\in V_i$ and $w\in V_j$ with $i<j$. 
\end{lemma}

The result of \citet{PTT99} generalises for upward 3D drawings as follows.

\begin{theorem}
\thmlabel{ColouredDrawing}
Every $n$-vertex $c$-colourable dag $G$ has an upward $c \times 4c^2n \times 4cn$ drawing with volume \Oh{c^4n^2}. 
\end{theorem}

\begin{proof}
Let $p$ be a prime number with $2c-1\leq p<4c$. Let $\{V_i:0\leq i\leq c-1\}$ be a $c$-colouring of $G$. Let $(v_1,v_2,\dots,v_n)$ be a topological ordering of $G$. Position each vertex $v_j\in V_i$ at $(\x_j,\y_j,\z_j)\in\mathbb{Z}^3$, where $\x_j:=i$ and $\y_j:=i\cdot\z_j$. It remains to compute the $\z_j$. If $v_1\in V_i$, then set $\z_1:=i^2\bmod{p}$. Now for all $j=2,3,\dots,n$, let $\z_j$ be the integer in $\{\z_{j-1}+1,\z_{j-1}+2,\dots,\z_{j-1}+p\}$ such that $\z_j\equiv i^2\pmod{p}$. Thus $\z_{j-1}<\z_j$. Hence arcs are upward, and no two vertices are mapped to the same point. By \lemref{Pach} with the \y- and \z-coordinates switched, the drawing is crossing-free. Since $0\leq \x_j\leq c-1$, the width is $c$. Since $\z_j\leq \z_{j-1}+p$, the height is at most $pn\leq4cn$. Since $\y_j<c\cdot\z_j$, the depth is less than $cpn\leq 4c^2n$.
\end{proof}

Many dags are $c$-colourable, for some constant $c$. These include dags whose
underlying undirected graph is outerplanar, is series parallel, is planar, or more generally, is from a proper minor-closed class, or has bounded degeneracy. \thmref{ColouredDrawing} implies that all such dags have upward 3D drawings with \Oh{n^2} volume.

\citet{PTT99} proved that the compete bipartite graph $K_{n,n}$ requires $\Omega(n^2)$ volume in every 3D drawing. Thus every acyclic orientation of
$K_{n,n}$ requires $\Omega(n^2)$ volume in every upward 3D drawing. Hence 
\thmref{ColouredDrawing} is tight for constant $c$. More generally, \citet{BCMW-JGAA04} proved that every 3D drawing of every $n$-vertex $m$-edge graph has at least $\frac{1}{8}(n+m)$ volume. 

%%%%%%%%%%%%%%%%%%%%%%%%%%%%%%%%%%%%%%%%%%%%%%%%%%%%%%%%%%%%%%%%%%%%%%%%%%%%%%
\subsection{Long Paths}
%%%%%%%%%%%%%%%%%%%%%%%%%%%%%%%%%%%%%%%%%%%%%%%%%%%%%%%%%%%%%%%%%%%%%%%%%%%%%%

We have the following lower bound, since every vertex in a directed path must be assigned a distinct \z-coordinate in an upward 3D drawing. 

\begin{lemma} 
\lemlabel{LongPathLowerBound}
Let $G$ be a dag that contains a directed path on $\ell$ vertices. Then the height of every upward 3D drawing of $G$ is at least $\ell$.\qed
\end{lemma}

Conversely, we have the following upper bound. 

\begin{theorem}
\thmlabel{LongPaths}
Every $n$-vertex dag $G$ with no directed path on $\ell$ vertices, has an upward $\Oh{\ell n}\times\Oh{n}\times\Oh{\ell}$ drawing with \Oh{\ell^2n^2} volume.
\end{theorem}

\begin{proof}
Colour each vertex $v\in V(G)$ by the number of vertices in the longest directed path ending at $v$. (This is well defined since $G$ is a dag.) The number of colours is at most $\ell$. Consider an arc $\arc{vw}\in A(G)$ such that $v$ is coloured $i$. Thus there is an $i$-vertex path $P$ ending at $v$. Moreover, $w\not\in P$ as otherwise $G$ would contain a directed cycle. Hence $(P,\arc{vw})$ is an $(i+1)$-vertex path ending at $w$. Thus the colour of $w$ is at least $i+1$. In particular, we have a proper $\ell$-colouring of $G$. 
The result follows from \lemref{PachPach} with the \x- and \z-coordinates switched.
\end{proof}

\thmref{LongPaths} is an improvement over \thmref{ColouredDrawing} whenever $\ell<\chi(G)^2$. 

%%%%%%%%%%%%%%%%%%%%%%%%%%%%%%%%%%%%%%%%%%%%%%%%%%%%%%%%%%%%%%%%%%%%%%%%%%%%%%
\section{Upward Track Layouts}
\seclabel{Track}
%%%%%%%%%%%%%%%%%%%%%%%%%%%%%%%%%%%%%%%%%%%%%%%%%%%%%%%%%%%%%%%%%%%%%%%%%%%%%%

%Loosely speaking, Cohen\etal\cite{CELR-Algo96} allow three `free' dimensions, whereas Pach\etal\cite{PTT99} use the assignment of vertices to colour classes to `fix' one dimension with two dimensions free. We use an assignment of vertices to tracks to fix two dimensions with one dimension free.  

Let $\{V_i:i\in I\}$ be a $t$-colouring of a graph $G$. Let $<_i$ be a total  order on each colour class $V_i$.  Then each pair $(V_i,<_i)$ is a \emph{track}, and $\{(V_i,<_i):i\in I\}$ is a $t$-\emph{track assignment} of $G$. To ease the notation we denote track assignments by $\{V_i:i\in I\}$  when the ordering on each colour class is implicit. An \emph{X-crossing} in a track assignment consists of two edges $vw$ and $xy$ such that $v<_ix$ and $y<_jw$, for distinct colours $i$ and $j$. A \emph{$t$-track layout} of $G$ is a $t$-track assignment of $G$ with no X-crossing.  The \emph{track-number} of $G$, denoted by \tn{G}, is the minimum integer $t$ such that $G$ has a $t$-track layout\footnote{Some authors \citep{DLMW-GD05, GLM-CGTA05, Giacomo-GD03, dGLW-CCCG02, DM-GD03} use a slightly different definition of track layout, in which \emph{intra-track} edges are allowed between consecutive vertices in a track. In keeping with the terminology of \citet{DMW-SJC05} and for consistency with the notion of an \emph{improper} colouring, we call this structure an \emph{improper} track layout, and use \emph{improper track-number} for the minimum number of tracks in this setting. The improper track-number is at most the track-number, and the track-number is at most twice the improper track-number \citep{DMW-SJC05}. Moreover, for every graph class $\mathcal{G}$ that includes all series parallel graphs, every graph in $\mathcal{G}$ has track-number at most some constant $t$ if and only if every graph in $\mathcal{G}$ has improper track-number at most $t$ \citep{DMW-SJC05}.}. 

Track layouts and track-number were introduced by \citet{DMW-SJC05} although they are implicit in many previous works \citep{FLW-JGAA03, HLR-SJDM92, HR-SJC92}. Track layouts and 3D drawings are closely related, as illustrated by the following results by Dujmovi{\'c}\etal\citep{DMW-SJC05, DujWoo-SubQuad-AMS}. 
\begin{theorem}[\citep{DMW-SJC05, DujWoo-SubQuad-AMS}]
\thmlabel{TrackDrawing}
Let $G$ be an $n$-vertex graph with chromatic number $\cn{G}\leq c$ and track-number $\tn{G}\leq t$. Then:
\begin{enumerate}
\item[\aaa] $G$ has a $\Oh{t}\times\Oh{t}\times\Oh{n}$ drawing with \Oh{t^2n} volume, and
\item[\bbb] $G$ has a $\Oh{c}\times\Oh{c^2t}\times\Oh{c^4n}$ drawing with \Oh{c^7tn} volume.
\end{enumerate}
Conversely, if a graph $G$ has an $X\times Y\times Z$ drawing, then $G$ has
track-number $\tn{G}\leq2XY$ \textup{(}and improper track-number at most $XY$\textup{)}.
\end{theorem}

The style of drawing produced by \thmref{TrackDrawing}(a) is illustrated in \figref{ThreeDimDrawing}.

\Figure{ThreeDimDrawing}{\includegraphics{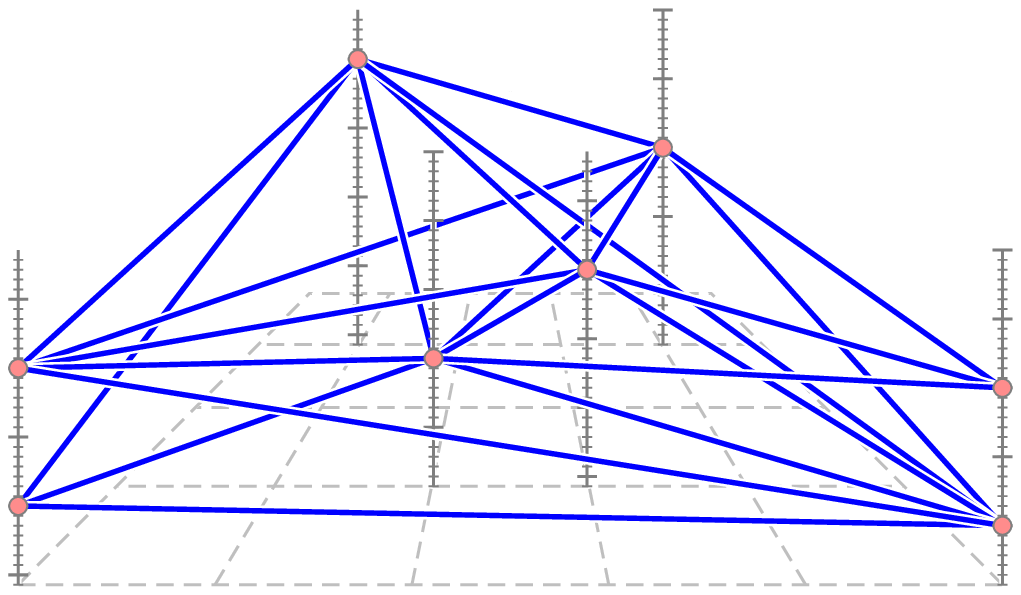}}{A 3D drawing produced from a $6$-track layout.}

The proof of \thmref{TrackDrawing}(a) implicitly used the following lemma.

\begin{lemma}[\citep{DMW-SJC05}]
\lemlabel{TrackDrawing}
Let $\{(V_i,<_i):1\leq i\leq t\}$ be a $t$-track layout of a graph $G$. Let $p>t$ be a prime number. Let $d_v$ be an integer for each vertex $v\in V(G)$, such that $d_v<d_w$ for all vertices $v,w\in V_i$ with $v<_iw$. If each vertex $v$ is placed at $(i,i^2\bmod{p},p\cdot d_v+i^3\bmod{p})$, then we obtain a \textup{(}crossing-free\textup{)} 3D drawing of $G$.
\end{lemma}

\citet{DLMW-GD05} extended the definition of track layouts to dags as follows\footnote{\citet{DLMW-GD05} allow intra-track arcs in their definition of upward track layout.}. An \emph{upward track layout} of a dag $G$ is a track layout of the underlying undirected graph of $G$, such that if $G^+$ is the directed graph obtained from $G$ by adding an arc from each vertex $v$ to the successor vertex in the track that contains $v$ (if it exists), then $G^+$ is still acyclic. The \emph{upward track-number} of $G$, denoted by \utn{G}, is the minimum integer $t$ such that $G$ has an upward $t$-track layout. \citet{DLMW-GD05} proved the following analogue of \thmref{TrackDrawing}(a).

\begin{theorem}[\citep{DLMW-GD05}]
\thmlabel{UpwardTrackDrawing}
Let $G$ be an $n$-vertex graph with upward track-number $\utn{G}\leq t$. Then $G$ has an upward $\Oh{t}\times\Oh{t}\times\Oh{tn}$ drawing with \Oh{t^3n} volume. Conversely, if a dag $G$ has an upward $X\times Y\times Z$ drawing then $G$ has upward track-number $\utn{G}\leq2XY$ \textup{(}and improper upward track-number at most $XY$\textup{)}.
\end{theorem}

\begin{proof}
Let $p$ be a prime number with $t<p\leq 2p$. For each vertex $v\in V(G)$, let $d_v$ be the maximum number of vertices in a directed path in $G^+$ that ends at $v$. Since each track induces a directed path in $G^+$, we have $d_v<d_w$ for all vertices $v$ and $w$ with $v<w$ in a single track. For each vertex $v$ in the $i$-th track, position $v$ at $(i,i^2\bmod{p},p\cdot d_v+i^3\bmod{p})$. Draw the arcs straight. By \lemref{TrackDrawing}, we obtain a crossing-free drawing. As in \thmref{LongPaths}, $d_v<d_w$ for every arc $\arc{vw}\in A(G)$. Thus the drawing is upward. The claimed volume bound holds since $d_v\leq n$. The converse results are proved in the same way as the converse results in \thmref{TrackDrawing}.
\end{proof}

For small values of $t$, the constants in \thmref{UpwardTrackDrawing} can be greatly improved.

\begin{lemma}
\lemlabel{ThreeTrackDrawing}
Every $n$-vertex dag $G$ that has an upward $3$-track layout $\{V_1,V_2,V_3\}$ has an upward $2\times 2\times n$ drawing with $4n$ volume.
\end{lemma}

\begin{proof}
Put the $i$-th vertex $v$ in a topological ordering of $G^+$ at $(0,0,i)$ if $v\in V_1$, at $(1,0,i)$ if $v\in V_2$, and at $(0,1,i)$ if $v\in V_3$. Draw each arc straight. Clearly we obtain an upward crossing-free drawing of $G$.
\end{proof}

It easily seen that \lemref{ThreeTrackDrawing} generalises to produce an upward $3\times 3\times n$ drawing of a $4$-track dag. This volume bound of $9n$ can be improved to $8n$ using a similar construction to one due to  \citet[Lemma~36]{DujWoo-DMTCS05}.

\begin{lemma}
\lemlabel{FourTrackDrawing}
Every $n$-vertex dag $G$ that has an upward $4$-track layout $\{V_1,V_2,V_3,V_4\}$ has an upward $2\times 2\times 2n$ drawing with $8n$ volume.
\end{lemma}

\begin{proof} Put the $i$-th vertex $v$ in a topological ordering of $G^+$ at $(0,0,2i)$ if $v\in V_1$,  at $(1,0,2i)$ if $v\in V_2$,  at $(0,1,2i)$ if $v\in V_3$,  and at $(1,1,2i-1)$ if $v\in V_4$. Draw each arc straight. Every arc is
upward, and the bounding box is at most $2\times2\times2n$. Suppose that edges $vw$ and $pq$ cross. Since there is no X-crossing in the track layout, $vw$ and $pq$ do not run between the same pair of tracks. The projection of the drawing onto the \x\y-plane is a subgraph of $K_4$ drawn with one crossing at $(\half,\half)$. The crossing is between the pairs of tracks $V_1V_4$ and $V_2V_3$. Thus without loss of generality $v\in V_1$, $w\in V_4$, $p\in V_2$, and $q\in V_3$. Hence the crossing point is $(\half,\half,\half(\z(v)+\z(w)))=(\half,\half,\half(\z(p)+\z(q)))$. This is a contradiction since $\z(v)+\z(w)$ is odd and $\z(p)+\z(q)$ is even. Hence the drawing is crossing-free. \end{proof}

%%%%%%%%%%%%%%%%%%%%%%%%%%%%%%%

\begin{lemma}
\lemlabel{FiveTrackDrawing}
Every $n$-vertex dag $G$ that has an upward $5$-track layout $\{V_1,V_2,V_3,V_4,V_5\}$ has an upward $4\times 4\times \frac{7}{5}n$ drawing with volume $22.4n$.
\end{lemma}

\begin{proof}
Without loss of generality, $|V_3|+|V_5|\leq\frac{2}{5}n$. Clearly we can assign distinct $\z$-coordinates to the vertices such that 
\begin{itemize}
\item $\z(v)<\z(w)$ for every arc $\arc{vw}\in A(G)$, 
\item $\z(v)$ is odd for every vertex $v\in V_3$, 
\item $\z(v)$ is even for every vertex $v\in V_5$, and 
\item $1\leq \z(v)\leq n+|V_3|+|V_5|\leq\frac{7}{5}n$ for every vertex $v\in V(G)$. 
\end{itemize}
As illustrated in \figref{UpwardTree}, put each vertex $v$ at $(1,1,\z(v))$ if $v\in V_1$, at $(2,3,\z(v))$ if $v\in V_2$, at $(2,4,\z(v))$ if $v\in V_3$, at $(3,2,\z(v))$ if $v\in V_4$, and at $(4,2,\z(v))$ if $v\in V_5$. Draw each arc straight. Every arc is upward, and the bounding box is at most $4\times4\times\frac{7}{5}n$. Suppose that edges $vw$ and $pq$ cross. Since there is no X-crossing in the track layout, $vw$ and $pq$ do not run between the same pair of tracks. The projection of the drawing onto the \x\y-plane is a subgraph of $K_5$ drawn with one crossing at $(\frac{8}{3},\frac{8}{3})$. The crossing is between the pairs of tracks $V_2V_5$ and $V_3V_4$. Thus without loss of generality $v\in V_2$, $w\in V_5$, $p\in V_3$, and $q\in V_4$. Now $vw$ and $pq$ intersect the line $\{(\frac{8}{3},\frac{8}{3},t):t\in\R\}$ respectively at $(\frac{8}{3},\frac{8}{3},\third(2\z(v)+\z(w)))$ and $(\frac{8}{3},\frac{8}{3},\third(2\z(q)+\z(p)))$. Thus $2\z(v)+\z(w)=2\z(q)+\z(p)$. This is a contradiction since $\z(w)$ is even and $\z(p)$ is odd. Hence the drawing is crossing-free.
\end{proof}

\Figure{UpwardTree}{\includegraphics{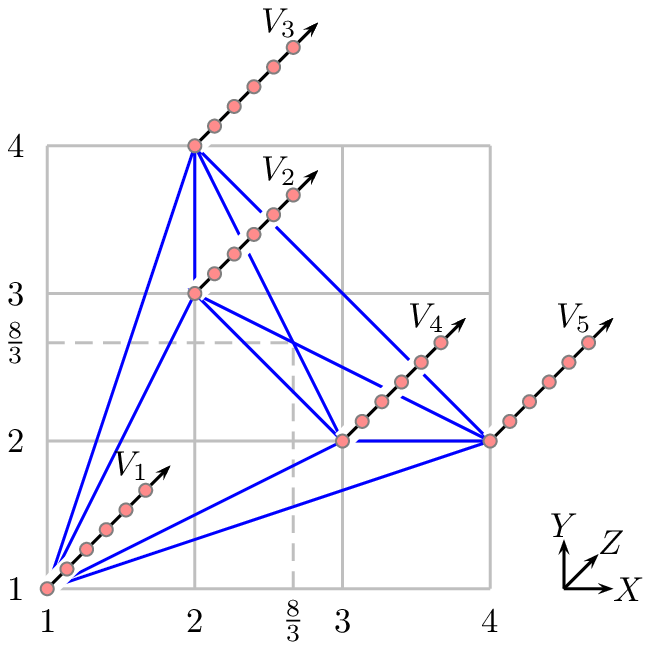}}{Construction of a 3D drawing from a $5$-track layout.}

%%%%%%%%%%%%%%%%%%%%%%%%%%%%%%%%%%%%%%%%%%%%%
\section{Strong Star Colourings} 
\seclabel{Strong}
%%%%%%%%%%%%%%%%%%%%%%%%%%%%%%%%%%%%%%%%%%%%%

\citet{DujWoo-SubQuad-AMS} defined a vertex colouring to be a \emph{strong star colouring} if between every pair of colour classes, all edges (if any) are incident to a single vertex. That is, each bichromatic subgraph consists of a star and possibly some isolated vertices. The \emph{strong star chromatic number} of a graph $G$, denoted by $\sscn{G}$, is the minimum number of colours in a strong star colouring of $G$. Note that \emph{star colourings}, in which each bichromatic subgraph is a star forest, have also been studied; see \cite{FRR-JGT04,NesOdM-03,Albertson-EJC04} for example. 

With an arbitrary order on each colour class in a strong star colouring, there is no X-crossing. Thus track-number $\tn{G}\leq\sscn{G}$, as observed by \citet{DujWoo-SubQuad-AMS}. Moreover, for a dag $G$ we can order each track by a topological vertex ordering of $G$, to obtain an upward track layout. Thus $\utn{G}\leq\sscn{G}$, as observed by \citet{DLMW-GD05}. 

\citet{DujWoo-SubQuad-AMS} proved that $\sscn{G}\leq 14\sqrt{\Delta m}$ and 
$\sscn{G}\leq 15m^{2/3}$ for every graph $G$ with maximum degree $\Delta$ and $m$ edges\footnote{Patrice Ossona de Mendez [personal communication] and 
J\'ean-Sebastien Sereni and St\'ephan Thomass\'e [personal communication] independently observed that if $H$ is the graph consisting of $k$ copies of $K_k$ (which has $m=\Theta(k^3)$ edges), then $\sscn{H}=\Theta(k^2)=\Theta(m^{2/3})$. Thus the general upper bound $\sscn{G}\leq\Oh{m^{2/3}}$ is best possible.}. In what follows we improve these bounds, by essentially replacing $\Delta$ by the weaker notion of degeneracy. The following concept will be useful. A colouring is \emph{harmonious} if every bichromatic subgraph has at most one edge; see \citep{Edwards97} for a survey on harmonious colourings. The \emph{harmonious chromatic number} of $G$, denoted by $h(G)$, is the minimum number of colours in a harmonious colouring of $G$. \citet{EM-JGT94} proved the following upper bound on $h(G)$. 

\begin{lemma}[\citep{EM-JGT94}]
\lemlabel{Harmonious}
Let $G$ be a $d$-degenerate graph with $m$ edges and maximum degree $\Delta$. Then $G$ has harmonious chromatic number $h(G)\leq 2\sqrt{2dm}+(2d-1)\Delta$.
\end{lemma}

\begin{lemma}
\lemlabel{SSCN}
Let $G$ be a $d$-degenerate graph $G$ with $m$ edges. Then the strong star chromatic number of $G$ satisfies $\sscn{G}\leq5\sqrt{2dm}$ and $\sscn{G}\leq(4+2\sqrt{2})m^{2/3}$.
\end{lemma}

\begin{proof}
For the first bound, let $A$ be the set of vertices of $G$ with degree at least $\sqrt{2m/d}$. Then $|A|\leq\sqrt{2dm}$. Now $G\setminus A$ has maximum degree at most $\sqrt{2m/d}$. Thus $h(G\setminus A)\leq2\sqrt{2dm}+(2d-1)\sqrt{2m/d}<
4\sqrt{2dm}$  by \lemref{Harmonious}. Using one colour for each vertex in $A$, we obtain a strong star colouring of $G$ with $5\sqrt{2dm}$ colours.

For the second bound, let $A$ be the set of vertices of $G$ with degree at least $m^{1/3}$. Then $|A|\leq 2m^{2/3}$. Now $G\setminus A$ has maximum degree and degeneracy at most $m^{1/3}$. Thus $h(G\setminus A)\leq2\sqrt{2m^{4/3}}+2m^{2/3}<(2+2\sqrt{2})m^{2/3}$ by \lemref{Harmonious}. Using one colour for each vertex in $A$, we obtain a strong star colouring of $G$ with $(4+2\sqrt{2})m^{2/3}$ colours.
\end{proof}

Since $\tn{G}\leq\sscn{G}$ and $\utn{G}\leq\sscn{G}$ we have the following corollary of \lemref{SSCN}.

\begin{corollary}
\corlabel{DegenTrackNumber}
Let $G$ be a $d$-degenerate graph with $n$ vertices and $m$ edges. 
Then the track-number of $G$ satisfies $\tn{G}\leq5\sqrt{2dm}<5d\sqrt{2n}$ and $\tn{G}\leq(4+2\sqrt{2})m^{2/3}$. The same bounds hold for the upward track-number of every acyclic orientation of $G$. \qed
\end{corollary}

\corref{DegenTrackNumber} does not give better bounds on the volume of upward 3D drawings than \thmref{ColouredDrawing} because of the cubic dependence on the upward track-number in \thmref{UpwardTrackDrawing}. However, for 3D drawings of undirected graphs, \thmref{TrackDrawing} and \corref{DegenTrackNumber} imply the following.

\begin{theorem}
\thmlabel{DegenVolume}
Let $G$ be an $n$-vertex graph with degeneracy $d$. Then $G$ has 
$\Oh{\sqrt{dm}}\times\Oh{\sqrt{dm}}\times\Oh{n}$ drawing with \Oh{dnm} volume. If $d$ is bounded, then $G$ has a $\Oh{1}\times\Oh{\sqrt{n}}\times\Oh{n}$ drawing with \Oh{n^{3/2}} volume.\qed
\end{theorem}

A number of notes on \thmref{DegenVolume} are in order.

\begin{itemize}

\item The above-mentioned bounds on the strong star chromatic number due to \citet{DujWoo-SubQuad-AMS} were non-constructive (since the proof was based on the Lov\'{a}sz Local Lemma). On the other hand, the proof of \lemref{Harmonious} is deterministic, and is easily seen to lead to polynomial time algorithms for  computing the colouring in \lemref{SSCN} and the drawing in 
\thmref{DegenVolume}.

\item The \Oh{dnm} volume bound in \thmref{DegenVolume} represents a qualitative improvement over the best previous comparable bound of \Oh{\Delta nm} in \citep{DujWoo-SubQuad-AMS}.

\item Graphs from a proper minor closed family were the largest class of graphs for which a \Oh{n^{3/2}} volume bound was previously known \citep{DujWoo-SubQuad-AMS}. The second part of \thmref{DegenVolume} is strictly stronger, since there are graph classes with bounded degeneracy but with unbounded clique minors. For example, the graph $K_n'$ obtained from $K_n$ by subdividing every edge once has degeneracy two, yet contains a $K_n$ minor. 

\end{itemize}

It is unknown what is the best possible bound on the track-number of graphs with bounded degeneracy. The graph $K_n'$ seems to be an important example. 

\begin{lemma}
$\tn{K_n'}=\Theta(n^{2/3})$.
\end{lemma}

\begin{proof}
First we prove the lower bound. Say $K_n'$ has a $t$-track layout. Some track contains a set $S$ of at least $p:=\ceil{n/t}$ `original' vertices of $K_n$. We can assume that $p$ is even. Say $S$ is ordered $v_1,v_2,\dots,v_p$ in this track. Let $T$ be the set of edges $v_iv_j$ of $K_n$ such that $1\leq i\leq p/2<j\leq p$. Observe that $|T|=p^2/4$. For all edges $e$ and $f$ in $T$, the division vertex of $K_n'$ that corresponds to $e$ and $f$ cannot be on the same track, as otherwise there will be an X-crossing. Thus the number of tracks $t\geq|T|\geq(n/t)^2/4$. Hence $t\geq(n/2)^{2/3}$.

Now we prove the upper bound. We can suppose that $p:=n^{1/3}$ is an integer. Partition the original vertices of $K_n'$ into $p^2$ sets, each with $p$ vertices. Places each set in its own track. Let $v_{i,k}$ be the $k$-th original vertex in the $i$-th track ($1\leq i\leq p^2$, $1\leq k\leq p$). For each such $k$, let $S_k:=\{v_{i,k}:1\leq i\leq p^2\}$. It remains to place the division vertices. As illustrated in \figref{CompletePrime}(a), the division vertices that correspond to edges with both endpoints in some $S_k$ are placed in one track ordered by increasing $k$. Since every vertex in each $S_k$ is in a distinct track, there is no X-crossing. For all $1\leq k<\ell\leq p$, let $T_{k,\ell}$ be the set $\{v_{i,k}v_{j,\ell}:1\leq i,j\leq p^2\}$ of edges of $K_n$. We place the division vertices of the edges in $T_{k,\ell}$ on two tracks as follows. Let $A_{k,\ell}:=\{v_{i,k}v_{j,\ell}\in T_{k,\ell}:1\leq i\leq j\leq p^2\}$ and $B_{k,\ell}:=\{v_{i,k}v_{j,\ell}\in T_{k,\ell}:1\leq j<i\leq p^2\}$. As illustrated in \figref{CompletePrime}(b), the division vertices that correspond to the edges in $A_{k,\ell}$ are placed on one track ordered by non-increasing $i$, breaking ties by decreasing $j$. The division vertices that correspond to the edges in $B_{k,\ell}$ are placed on one track ordered by non-decreasing $j$, breaking ties by increasing $i$. It is easily seen that there is no X-crossing. In total we have $p^2+1+2\binom{p}{2}\sim 2p^2\sim 2n^{2/3}$ tracks. (Note that the constant in this upper bound can be improved by combining tracks $A_{k,\ell}$ and $A_{k',\ell'}$ when $k<\ell<k'<\ell'$; we omit the details.)
\end{proof}

\Figure{CompletePrime}{\includegraphics{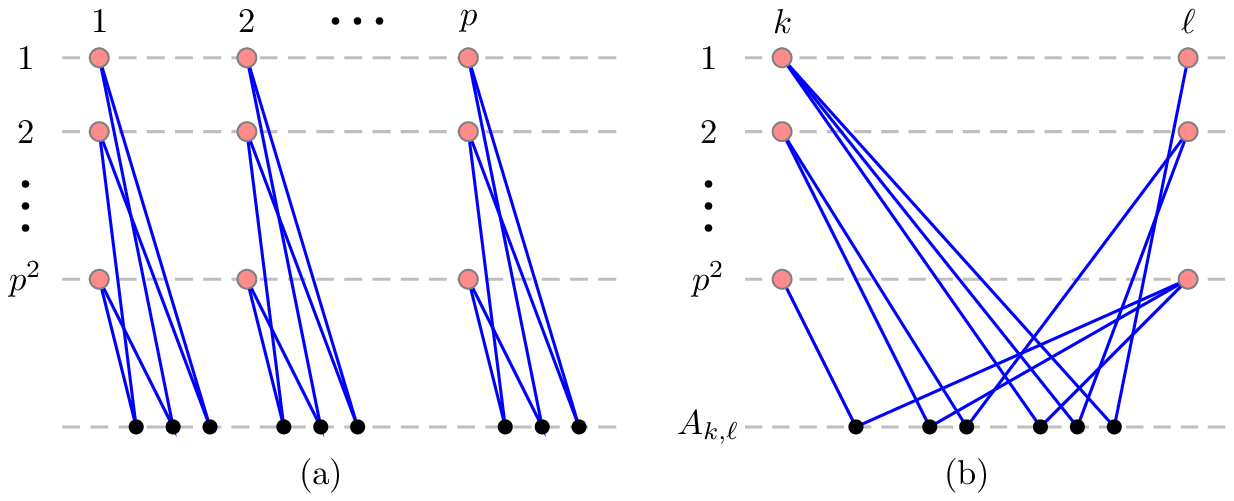}}{Construction of track layout of $K_n'$.}

\comment{more detail at the end of this proof, draw a figure}

%%%%%%%%%%%%%%%%%%%%%%%%%%%%%%%%%%%%%%%%%%%%%%%%%%%%%%%%%%%%%%%%%%%%%%%%%%%%%%
\section{Upward Queue Layouts}
\seclabel{Queue}
%%%%%%%%%%%%%%%%%%%%%%%%%%%%%%%%%%%%%%%%%%%%%%%%%%%%%%%%%%%%%%%%%%%%%%%%%%%%%%

A $k$-\emph{queue layout} of a graph $G$ consists of a vertex ordering $\sigma$ of $G$, and a partition $\{E_1,E_2,\dots,E_k\}$ of $E(G)$, such that no two edges in each $E_i$ are \emph{nested} in $\sigma$. That is, for all edges $vw,xy\in E_i$, we do not have $\sigma(v)<\sigma(x)<\sigma(y)<\sigma(w)$. The \emph{queue-number} of $G$, denoted by \qn{G}, is the minimum integer $k$ such that $G$ has a $k$-queue layout. Queue layouts and queue-number were introduced by Heath\etal\citep{HLR-SJDM92,HR-SJC92}; see \citep{DujWoo-DMTCS04} for references and results.

Heath\etal\citep{HP-SJC99, HPT-SJC99} extended the definition of queue layouts to dags as follows. A \emph{upward} $k$-\emph{queue layout} of a dag $G$ is a $k$-\emph{queue layout} of the underlying undirected graph of $G$ such that the vertex ordering $\sigma$ is topological. For example, every tree has a $1$-queue layout \citep{HR-SJC92}, and every tree dag has an upward $2$-queue layout \citep{HPT-SJC99}. The \emph{upward queue-number} of $G$, denoted by \uqn{G}, is the minimum integer $k$ such that $G$ has an upward $k$-queue layout. 

Consider a vertex colouring $\{V_i:i\in\Z\}$ of a graph $G$ in which the colours are integers. An edge $vw\in E(G)$ with $v\in V_i$ and $w\in V_j$ has \emph{span} $|j-i|$. This definition naturally extends to track layouts $\{(V_i,<_i):i\in\Z\}$. Let $(\dots,V_{-1},V_0,V_1,\dots)$ denote the vertex ordering $\sigma$ in which $\sigma(v)<\sigma(w)$ whenever $v\in V_i$ and $w\in V_j$ for some $i<j$, or $v<_i w$ within some $V_i$. \citet[Lemma~19]{DPW-DMTCS04} proved a characterisation of $1$-queue graphs in terms of track layouts with span two that immediately generalises for upward $1$-queue layouts as follows.  

\begin{lemma}
\lemlabel{UpwardOneQueue}
A dag $G$ has an upward $1$-queue layout if and only if $G$ has a track layout $\{V_i:i\in\Z\}$, such that for every arc $\arc{vw}\in A(G)$ with $v\in V_i$ and $w\in V_j$, we have $i<j\leq i+2$, and if $j=i+2$ then $w$ is the first vertex in $V_j$, and there is no arc \arc{xy} with $v<x\in V_i$ and $y\in V_{i+1}$. In particular, $(\dots,V_{-1},V_0,V_1,\dots)$ defines an upward $1$-queue layout of $G$.
\end{lemma}

An important technique for `wrapping' an undirected track layout is generalised for a particular type of upward track layout as follows. 

\begin{lemma}
\lemlabel{Wrap}
Let $\{V_i:i\in\Z\}$ be a track layout of a dag $G$, such that for every arc $\arc{vw}\in A(G)$ with $v\in V_i$ and $w\in V_j$, we have $i<j\leq i+s$.
\begin{enumerate}
\item[\aaa] Then $\{V_i:i\in\Z\}$ is an upward track layout.
\item[\bbb] The vertex ordering $\sigma=(\dots,V_{-1},V_0,V_1,\dots)$ defines an upward $s$-queue layout of $G$, and $G$ has upward queue-number $\uqn{G}\leq s$.
\item[\ccc] For each $0\leq i\leq 2s$, define 
\begin{equation*}
W_i:=(\dots,V_{i-2(2s+1)},V_{i-(2s+1)},V_i,V_{i+2s+1},V_{i+2(2s+1)},\dots).
\end{equation*} 
Then $\{W_0,W_1,\dots,W_{2s}\}$ is an upward $(2s+1)$-track layout of $G$, and $G$ has upward track-number $\utn{G}\leq2s+1$.
\end{enumerate}
\end{lemma} 
 
\begin{proof} 
Observe that $\sigma$ is a topological ordering of $G^+$. Thus $G^+$ is acyclic, and $\{V_i:i\in\Z\}$ is an upward track layout. This proves (a). Now we prove (b). Two arcs with the same span are not nested in $\sigma$ \citep[Lemma~5.2]{DMW-SJC05}. Thus we can partition the arcs into $s$ queues in $\sigma$ according to their span. Now we prove (c). The track assignment $\{W_0,W_1,\dots,W_{2s}\}$ is upward since the corresponding graph $G^+$ is acyclic. \citet[Lemma~3.4]{DMW-SJC05} proved that there is no X-crossing in $\{W_0,W_1,\dots,W_{2s}\}$. 
\end{proof}

A \emph{$(k,t)$-track layout} of a graph $G$ is a $t$-track assignment of $G$ in which every edge is assigned one of $k$ colours, such that there is no monochromatic X-crossing. By \tn[k]{G} we denote the minimum integer $t$ such that $G$ has a $(k,t)$-track layout. Thus $\tn{G}=\tn[1]{G}$. These definitions immediately generalise to the setting of upward $(k,t)$-track layouts.

\citet{DPW-DMTCS04} proved that queue-number and track-number are tied, in the sense that there is a function $f$ such that for every graph $G$, we have $\tn{G}\leq f(\qn{G})$ and $\qn{G}\leq f(\tn{G})$. In one direction the proof is easy. Given a $(k,t)$-track layout $\{V_1,V_2,\dots,V_t\}$ of $G$, \citet{DPW-DMTCS04} proved that the vertex ordering $(V_1,V_2,\dots,V_t)$ admits a $k(t-1)$-queue layout, and thus $\qn{G}\leq k(\tn[k]{G}-1)$. \citet{DLMW-GD05} proved the case $k=1$ of the following analogous relationship between upward track-number and upward queue-number.

\begin{lemma}
\lemlabel{Track2Queue}
For all $k\geq1$ and for every graph $G$, 
\begin{equation*}
\uqn{G}\leq k\cdot\binom{\utn[k]{G}}{2}.
\end{equation*}
\end{lemma}

\begin{proof}
Let $\sigma$ be a topological ordering of $G^+$ (defined with respect to a given $(k,t)$-track layout). Thus the order of each track is preserved in $\sigma$. Monochromatic arcs between each pair of tracks form a queue in $\sigma$ \citep[Lemma~14]{DPW-DMTCS04}.
\end{proof}

%The undirected version of the following converse result was proved by \citet{DPW-DMTCS04}.

%\begin{lemma} If a $c$-colourable graph $G$ has an upward $k$-queue layout, then $G$ has an upward $(2k,c)$-track layout. \end{lemma}

%\begin{proof} Let $\sigma$ be a topological ordering that admits an upward $k$-queue layout of $G$. Let $\{V_1,V_2,\dots,V_c\}$ be a $c$-colouring of $G$. Consider each $V_i$ to be a track ordered by $\sigma$. Then arcs \arc{vw} and \arc{xy} with $v,x\in V_i$ and $w,y\in V_j$ do not form an X-crossing. Similarly, arcs \arc{vw} and \arc{xy} with $v,y\in V_i$ and $w,x\in V_j$ do not form an X-crossing. Thus we obtain an upward $(2k,c)$-track layout of $G$. \end{proof}

\begin{lemma}
\lemlabel{Queue2Track}
Upward track-number is bounded by upward queue-number. In particular, every dag $G$ with upward queue-number $\uqn{G}\leq q$ has upward track-number $\utn{G}\leq4q\cdot4^{q(2q-1)(4q-1)}$.
\end{lemma}

\begin{proof} 
\citet[Theorem~8]{DPW-DMTCS04} proved that every (undirected) graph $G$ with queue-number $\qn{G}\leq q$ has track-number $\tn{G}\leq4q\cdot4^{q(2q-1)(4q-1)}$. In the proof, the ordering of $V(G)$ in the given $q$-queue layout is preserved in every track of the track layout. Thus the result also holds for upward track-number. 
\end{proof}

\twolemref{Track2Queue}{Queue2Track} imply 

\begin{theorem}
\thmlabel{Tied}
Upward track-number and upward queue-number are tied.
\end{theorem}

%%%%%%%%%%%%%%%%%%%%%%%%%%%%%%%%%%%%%%%%%%%%%%%%%%%%%%%%%%%%%%%%%%%%%%%%%%%%%%
\section{An Example}
\seclabel{Example}
%%%%%%%%%%%%%%%%%%%%%%%%%%%%%%%%%%%%%%%%%%%%%%%%%%%%%%%%%%%%%%%%%%%%%%%%%%%%%%

As illustrated in \figref{Example}, let $G_n$ be the dag with vertex set $\{u_i:1\leq i\leq 2n\}$ and arc set $\{\arc{u_iu_{i+1}}:1\leq i\leq 2n-1\}\cup\{\arc{u_iu_{2n-i+1}}:1\leq i\leq n\}$.

\Figure{Example}{\includegraphics{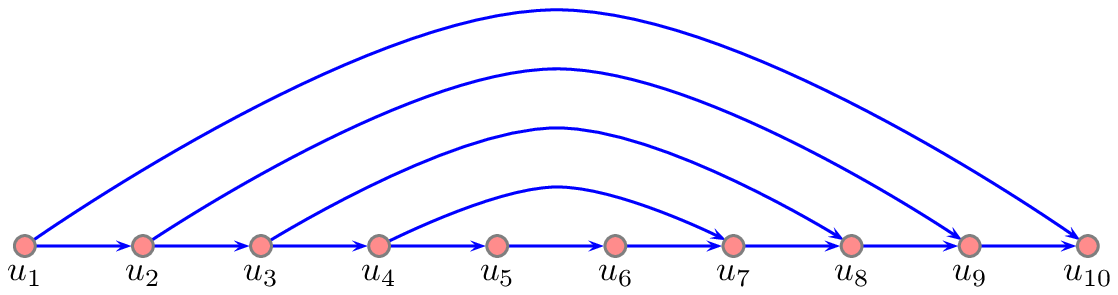}}{Illustration of $G_5$.}

Observe that $G_n$ is outerplanar and has a Hamiltonian directed path $(u_1,u_2,\dots,u_{2n})$. Thus $(u_1,u_2,\dots,u_{2n})$ is the only topological ordering of $G_n$, in which the edges $\{u_iu_{2n-i+1}:1\leq i\leq n\}$ are pairwise nested. Thus $\uqn{G_n}\geq n$; it is easily seen that in fact $\uqn{G_n}=n$. These observation were made by \citet{HPT-SJC99}. 

\thmref{Tied} implies that $G_n$ has unbounded upward track-number. \citet{DLMW-GD05} proved the same result with the much stronger bound of $\utn{G_n}\geq\sqrt{2n}$, which follows from \lemref{Track2Queue} with $k=1$ and since $\uqn{G_n}\geq n$. An upper bound of $\utn{G_n}\leq\Oh{\sqrt{n}}$ follows from \corref{DegenTrackNumber} (and since $G_n$ has bounded degree, from the earlier bounds on track-number in \citep{DujWoo-SubQuad-AMS}). \citet{DLMW-GD05} gave an elegant construction of an improper \Oh{\sqrt{n}}-track layout of $G_n$. 
Suppose that $G_n$ has an upward $X\times Y\times Z$ drawing. Then
$Z\geq2n$ by \lemref{LongPathLowerBound}. The second part of \thmref{UpwardTrackDrawing} implies that $2XY\geq\utn{G_n}\geq\sqrt{2n}$. Hence the volume is $\Omega(n^{3/2})$, as proved by \citet{DLMW-GD05}. This result highlights a substantial difference between 3D drawings of undirected graphs and upward 3D drawings of dags, since every (undirected) outerplanar graph has a 3D drawing with linear volume \citep{FLW-JGAA03}. In the full version of their paper, \citet{DLMW-GD05} constructed an upward 3D drawing of $G_n$ with \Oh{n^{3/2}} volume. It is unknown whether every $n$-vertex outerplanar dag has 
an upward 3D drawing with \Oh{n^{3/2}} volume. 

%\comment{I suspect there is an upward $2\times\Oh{\sqrt{n}}\times\Oh{n}$ drawing of $G_n$ with \Oh{n^{3/2}} volume. Assume that $r:=\sqrt{n}$ is an integer. Rename the vertices of $G_n$ as follows. For all $i$ and $j$ with $0\leq i,j\leq r-1$, let $v_{i,j}:=u_{ir+j+1}$ and  $w_{i,j}:=u_{r^2+ir+j+1}$. Observe that the arcs not in the Hamiltonian path of $G_n$ are of the form \arc{v_{i,j}w_{r-1-i,r-1-j}}. Place each vertex $v_{i,j}$ at $(0,i,ri+j)$, and place each vertex $w_{i,j}$ at $(1,r-j,r^2+ir+j)$ or something like this. }

%%%%%%%%%%%%%%%%%%%%%%%%%%%%%%%%%%%%%%%%%%%%%%%%%%%%%%%%%%%%%%%%%%%%%%%%%%%%%%
\section{Upward Layouts of Trees}
\seclabel{Trees}
%%%%%%%%%%%%%%%%%%%%%%%%%%%%%%%%%%%%%%%%%%%%%%%%%%%%%%%%%%%%%%%%%%%%%%%%%%%%%%

A \emph{caterpillar} is a tree, such that deleting the leaves gives a path (called the \emph{spine}). A graph has a $2$-track layout if and only if it is a caterpillar forest \citep{HS72}. \citet{FLW-JGAA03} proved that every tree has a $3$-track layout, which is best possible for every tree that is not a caterpillar. Using this result, \citet{DLMW-GD05} proved that every tree dag has an upward $7$-track layout, and that there exist tree dags with no upward $3$-track layout. We improve this upper bound as follows.

\comment{Characterise the dags with an upward 2-track layouts}

\begin{theorem}
\thmlabel{TreeTrackLayout}
Every tree dag has an upward $5$-track layout.
\end{theorem}

\thmref{TreeTrackLayout} follows from \lemref{Wrap}(c) with $s=2$ and the following lemma. Similarly, \lemref{Wrap}(b) with $s=2$ and following lemma imply that every tree dag has upward queue-number at most two, as proved by \citet{HPT-SJC99}.

\begin{lemma} 
\lemlabel{TreeUpwardTrack}
Every tree dag $T$ has a track layout $\{V_i:i\in\Z\}$, such that for every arc $\arc{vw}\in A(T)$ with $v\in V_i$ and $w\in V_j$, we have $i<j\leq i+2$. 
\end{lemma}

\begin{proof} 
Choose an arbitrary vertex $r$ of $T$. Consider a vertex $v$. The \emph{distance} $d(v)$ is the distance between $v$ and $r$ in the underlying undirected tree of $T$. Let $a(v)$ be the number of arcs on the path from $v$ to $r$ that are directed toward $r$. Let $b(v)$ be the number of arcs on the path from $v$ to $r$ that are directed away from $r$. Note that $d(v)=a(v)+b(v)$. As illustrated in \figref{TreeTrack}, put $v$ in track $V_{2b(v)-a(v)}$. Within each track the vertices are ordered in non-decreasing order of distance from $r$. It remains to order the vertices in a single track at the same distance. We do so inductively by increasing distance. Suppose we have ordered all the vertices at distance at most $d-1$. Consider two vertices $v_1,v_2$ in the same track and at distance $d$. Let $w_1,w_2$ be their respective parent vertices at distance $d-1$. If $w_1<w_2$ in the same track, then place $v_1<v_2$. If $w_1$ and $w_2$ are in different tracks, then the relative order of $v_1$ and $v_2$ is not important. This completes the construction. 

\Figure{TreeTrack}{\includegraphics{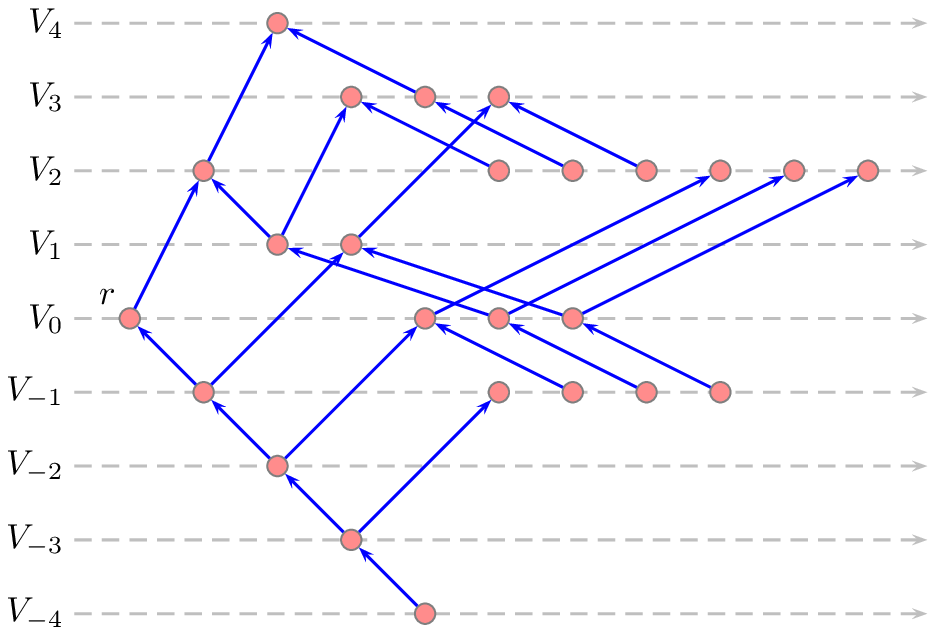}}{An upward track layout of a tree dag with span two.}

Consider an arc $\arc{vw}\in A(T)$ with $v\in V_i$ and $w\in V_j$. First suppose that $\arc{vw}$ is directed toward $r$. Then $b(w)=b(v)$ and $a(v)=a(w)+1$. Thus $j=i+1$ and $\arc{vw}$ has span one. Now suppose that $\arc{vw}$ is directed away from $r$. Then $a(w)=a(v)$ and $b(w)=b(v)+1$. Thus $j=i+2$ and $\arc{vw}$ has span two. It remains to prove that there is no X-crossing. Consider arcs \arc{vw} and \arc{xy} between the same pair of tracks. Thus both arcs have the same span. First suppose their span is one. Thus both \arc{vw} and \arc{xy} are directed towards $r$. Hence $d(w)=d(v)+1$ and $d(y)=d(x)+1$. Without loss of generality $v<x$ in their track. By construction $d(v)\leq d(w)$. If $d(v)<d(x)$ then $d(w)<d(y)$, and $w<y$ in their track since tracks are ordered by non-decreasing distance. If $d(v)=d(x)$ then $d(w)=d(y)$, and by construction $w<y$ in their track. In both cases, the arcs do not form an X-crossing. If \arc{vw} and \arc{xy} have span two, then both \arc{vw} and \arc{xy} are directed away from $r$, and an analogous argument proves that the arcs do not form an X-crossing. 
\end{proof}

\citet{DLMW-GD05} proved that every tree dag has an upward $7\times7\times7n$ drawing. \lemref{FiveTrackDrawing} and \thmref{TreeTrackLayout} imply the following improved bound.

\begin{theorem}
Every $n$-vertex tree dag has an upward $4\times4\times\frac{7}{5}n$ drawing.\qed
\end{theorem}

The next result generalises and improves upon the result of \citet{DLMW-GD05} that every directed path has an improper upward $3$-track layout.

\begin{theorem} 
\thmlabel{Caterpillar}
Every caterpillar dag has an upward $3$-track layout and an upward $1$-queue layout.
\end{theorem}

\thmref{Caterpillar} follows from the following lemma and \lemref{Wrap} with $s=1$.

\begin{lemma}
Every caterpillar dag $C$ has a track layout $\{V_i:i\in\Z\}$ such that for every arc $\arc{vw}\in A(C)$, if $v\in V_i$ then $w\in V_{i+1}$. 
\end{lemma}

\begin{proof}
Let $P=(v_1,v_2,\dots,v_n)$ be the spine of $C$. Put $v_1$ in track $V_0$. For $j=2,3,\dots,n$, place $v_j$ as follows. Say $v_{j-1}\in V_i$. If $v_{j-1}v_j\in A(P)$ then put $v_j$ rightmost in track $V_{i+1}$. Otherwise $v_jv_{j-1}\in A(P)$, in which case put $v_j$ rightmost in track $V_{i-1}$. Clearly $\{V_i:i\in\Z\}$ is a track layout of $P$. Consider a leaf vertex $w$ adjacent to $v_j\in V_i$. If $\arc{wv_j}\in A(C)$ then put $w$ in $V_{i-1}$. Otherwise $\arc{v_jw}\in A(C)$, in which case put $w$ in $V_{i+1}$. Order each track as shown in \figref{Caterpillar}. Clearly $\{V_i:i\in\Z\}$ is a track layout of $C$. 
\end{proof}

\Figure{Caterpillar}{\includegraphics{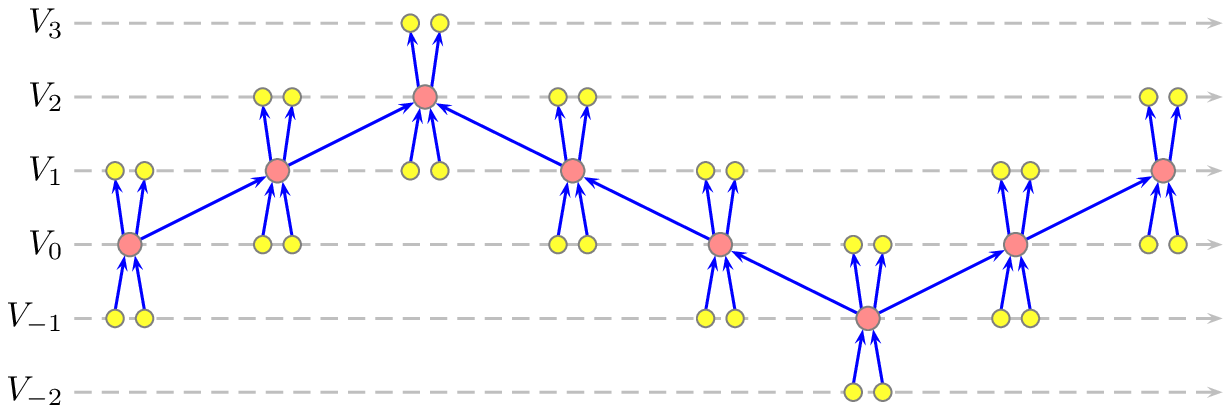}}{An upward track layout of a caterpillar dag with span one.}

Note the converse result.

\begin{lemma}
Suppose that every orientation of a tree $T$ has an upward track layout $\{V_i:i\in\Z\}$ such that for every arc $\arc{vw}$, if $v\in V_i$ then $w\in V_{i+1}$. Then $T$ is a caterpillar.
\end{lemma}

\begin{proof}
On the contrary, suppose that $T$ is a tree that is not a caterpillar, and every orientation of $T$ has the desired track layout. Then $T$ contains a \emph{2-claw} \citep{HS72}, which is the tree with vertices $r,u,v,w,x,y,z$ and edges $ru,rv,rw,ux,vy,wz$. Orient the edges \arc{ru}, \arc{rv}, \arc{rw}, \arc{xu}, \arc{yv}, \arc{zw}. Then this directed $2$-claw has the desired upward track layout. Say $r\in V_i$. Then $u,v,w\in V_{i+1}$ and $x,y,z\in V_i$. Without loss of generality, $v$ is between $u$ and $w$ in $V_{i+1}$. Thus \arc{yv} forms an X-crossing with either \arc{xu} or \arc{zw}, which is the desired contradiction.
\end{proof}

\lemref{ThreeTrackDrawing} and \thmref{Caterpillar} imply the following.

\begin{corollary} 
Every $n$-vertex caterpillar dag has a $2\times 2\times n$ upward 3D drawing.\qed
\end{corollary}

%%%%%%%%%%%%%%%%%%%%%%%%%%%%%%%%%%%%%%%%%%%%%%%%%%%%%%%%%%%%%%%%%%%%%%%%%%%%%%
\section{Upward Layouts of Subdivisions}
\seclabel{Subdivisions}
%%%%%%%%%%%%%%%%%%%%%%%%%%%%%%%%%%%%%%%%%%%%%%%%%%%%%%%%%%%%%%%%%%%%%%%%%%%%%%

A \emph{subdivision} of a graph $G$ is a graph obtained from $G$ by replacing each edge $vw\in E(G)$ by a path with endpoints $v$ and $w$. Each vertex of $G$ is called an \emph{original} vertex of the subdivision. \citet{DujWoo-DMTCS05} proved that every graph $G$ has a $2$-queue subdivision with \Oh{\log\qn{G}} division vertices per edge, and that this bound is best possible. We now prove a similar result for upward queue layouts. A \emph{subdivision} of a dag $G$ is a graph obtained from $G$ by replacing each arc $\arc{vw}\in A(G)$ by a path from $v$ to $w$. 

Let $\sigma$ be a vertex ordering of a graph $G$. The \emph{bandwidth} of $\sigma$ is 
\begin{equation*}
\max\{|\sigma(w)-\sigma(v)|:vw\in E(G)\}.
\end{equation*} 
The \emph{bandwidth} of $G$ is the minimum bandwidth of a vertex ordering of $G$. The \emph{directed bandwidth} of a dag $G$ is the minimum bandwidth of a topological vertex ordering of $G$. 

\begin{theorem}
\thmlabel{TwoQueueSub}
Every dag $G$ with directed bandwidth $b$ has an upward $2$-queue subdivision with at most $\half(b-1)$ division vertices per edge.
\end{theorem}

\begin{proof}
Let $\sigma=(v_1,v_2,\dots,v_n)$ be a topological ordering of $G$ with bandwidth $b$. For every arc $\arc{v_iv_j}\in A(G)$ with $j-i$ even, replace \arc{v_iv_j} by the directed path 
\begin{equation*}
(v_i,x(i,j,i+2),x(i,j,i+4),\dots,x(i,j,j-4),x(i,j,j-2),v_j)
\end{equation*} 
with $\half(j-i-2)$ division vertices. For every arc $\arc{v_iv_j}\in A(G)$ with $j-i$ odd, replace \arc{v_iv_j} by the directed path 
\begin{equation*}
(v_i,x(i,j,i+1),x(i,j,i+3),\dots,x(i,j,j-4),x(i,j,j-2),v_j)
\end{equation*} 
with $\half(j-i-1)$ division vertices. 

First we construct an upward $(2,n)$-track layout of this subdivision of $G$. Place each vertex $v_i$ leftmost in track $V_i$. Position every vertex $x(i,j,\ell)$ in track $V_\ell$. It remains to order the vertices within each track. Observe that $V_1=\{v_1\}$. Order $V_2$ arbitrarily (with $v_2$ leftmost). Now order $V_3$, then $V_4$, and so on up to $V_n=\{v_n\}$. In track $V_\ell$, except for the original vertex $v_\ell$, each vertex has exactly one neighbour in $\cup\{V_i:1\leq i<\ell\}$. Thus we can order $V_\ell$ with $v_\ell$ leftmost, so that no two edges incident to vertices in $V_\ell\setminus\{v_\ell\}$ form an X-crossing. The only possible X-crossing involves an edge incident to $v_\ell$. Thus we can colour all incoming edges $x(i,\ell,\ell-2)v_\ell$ \emph{green}, colour all other edges \emph{blue}, and there is no monochromatic X-crossing. Hence we have a $(2,n)$-track layout.

It is easily verified that the blue edges satisfy \lemref{UpwardOneQueue}  (with respect to the upward $(2,n)$-track layout $\{V_1,V_2,\dots,V_n\}$), as do the green edges. Thus the vertex ordering $(V_1,V_2,\dots,V_n)$ defines an upward $2$-queue layout of $G$.
\end{proof}

\Figure{TwoTrack}{\includegraphics{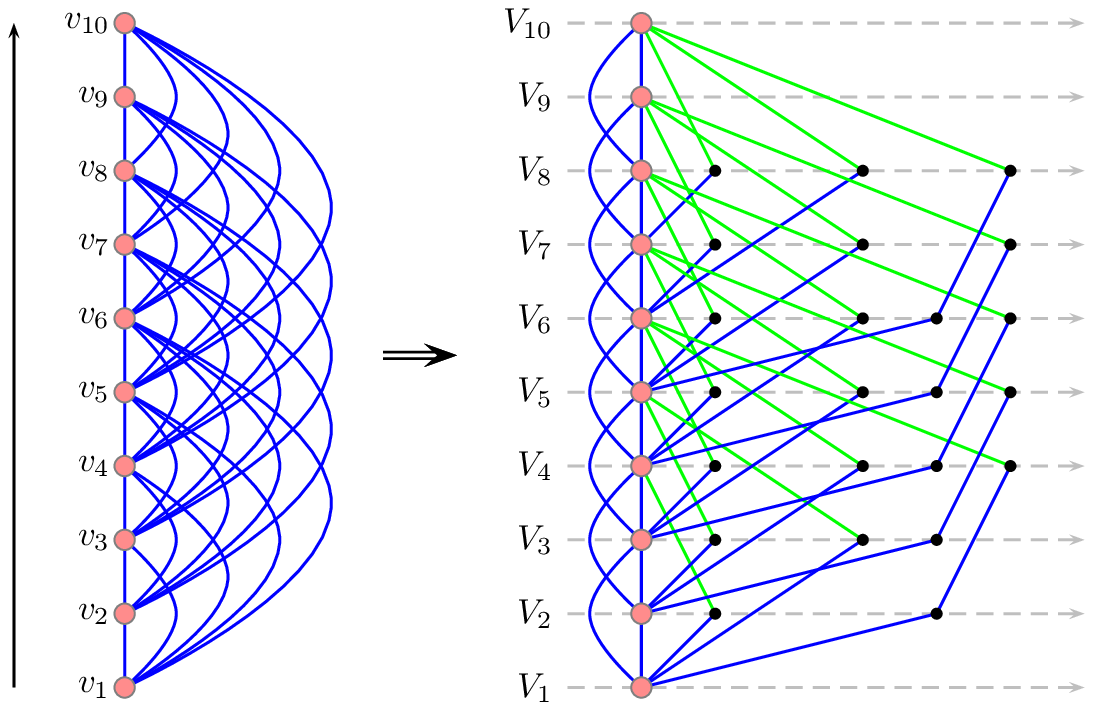}}{$(2,n)$-track layout in the proof of \thmref{TwoQueueSub}.}

\citet{DujWoo-DMTCS05} proved that every graph has a $4$-track subdivision with $\Oh{\log\qn{G}}$ division vertices per edge, and that this bound is best possible. We have the following similar result for upward track layouts. 

%\begin{theorem}\thmlabel{FourTrackSub}Every dag $G$ with directed bandwidth $b$ has an upward $4$-track subdivision with at most $b$ division vertices per edge.\end{theorem}

%\begin{proof}Let $\sigma=(v_1,v_2,\dots,v_n)$ be a topological ordering of $G$ with bandwidth $b$. Replace every arc $\arc{v_iv_j}\in A(G)$ by the directed path \begin{equation*}(v_i,x(i,j,i+1),x(i,j,i+2),\dots,x(i,j,j),v_j),\end{equation*} which has at most $b$ division vertices. Let $W=\{x(i,j,\ell):v_iv_j\in A(G),\ell=j\}$. Observe that $W$ is independent. Consider $W$ to be a track ordered firstly by non-decreasing $\ell$, breaking ties by increasing $i$. Let $x(i,j,i)$ denote the vertex $v_i$ for all arcs $\arc{v_iv_j}\in A(G)$. For each $m\in\{0,1,2\}$, define \begin{equation*}V_m:=\{x(i,j,\ell):\ell\equiv m\pmod{3},v_iv_j\in A(G),i\leq\ell\leq j-1\}.\end{equation*} Observe that $V_m$ is independent. Consider $V_m$ to be a track ordered firstly by non-decreasing $\ell$, breaking ties by non-increasing $i$, and then by by increasing $j$. Clearly every vertex of the subdivision is in one of the four tracks, and the track assignment is upward. A straightforward case-analysis proves that there is no X-crossing; cf.~\lemref{Wrap}(c). \end{proof}

%%%%%%%%%%%%%%%%%%%%%%%%%%%%%%%%%%%%%%%%%%%

\begin{theorem}
\thmlabel{FourTrackSub}
Every dag $G$ with directed bandwidth $b$ has an upward $4$-track subdivision with at most $b$ division vertices per arc.
\end{theorem}

\begin{proof} 
The proof proceeds in three steps. First, we construct an upward $(2,n)$-track layout of a subdivison $G'$ of $G$ that has at most $b-1$ division vertices per arc. Second, we construct an upward $(n+1)$-track layout of a subdivison $G''$ of $G$ that has at most $b$ division vertices per arc. Finally, we wrap $n$ of the tracks into three tracks, to obtain an upward $4$-track layout of $G''$.

Let $(v_0,v_1,\dots,v_{n-1})$ be a topological ordering of $G$ with bandwidth $b$. Let $G'$ be the subdivision of $G$ obtained by replacing each arc $\arc{v_iv_j}\in A(G)$ for which $j\geq i+2$ by the directed path 
\begin{equation*}
(v_i, x(i, j, i + 1), x(i, j, i + 2), \dots, x(i, j, j-1), v_j).
\end{equation*}
Note that $G'$ has at most $b-1$ division vertices per arc. Let $g(i,j)$ denote the final arc $(x(i,j,i-1),v_j)$ in each such path. Colour $g(i,j)$ \emph{green}, and colour the remaining arcs in $G'$ \emph{blue} (including the non-subdivided arcs $\arc{v_iv_{i+1}}$). 

As illustrated in \figref{FourTrack}(b), create a $(2,n)$-track layout of $G'$ as follows. Place each original vertex $v_\ell$ leftmost in track $V_\ell$ followed by the division vertices $x(i,j,\ell)$ for all arcs $\arc{v_iv_j}\in A(G)$ with $i<\ell<j$. Order the division vertices in track $V_\ell$ by non-increasing $i$, breaking ties by increasing $j$. It is simple to verify that no two monochromatic edges form an X-crossing. 

\Figure{FourTrack}{\includegraphics{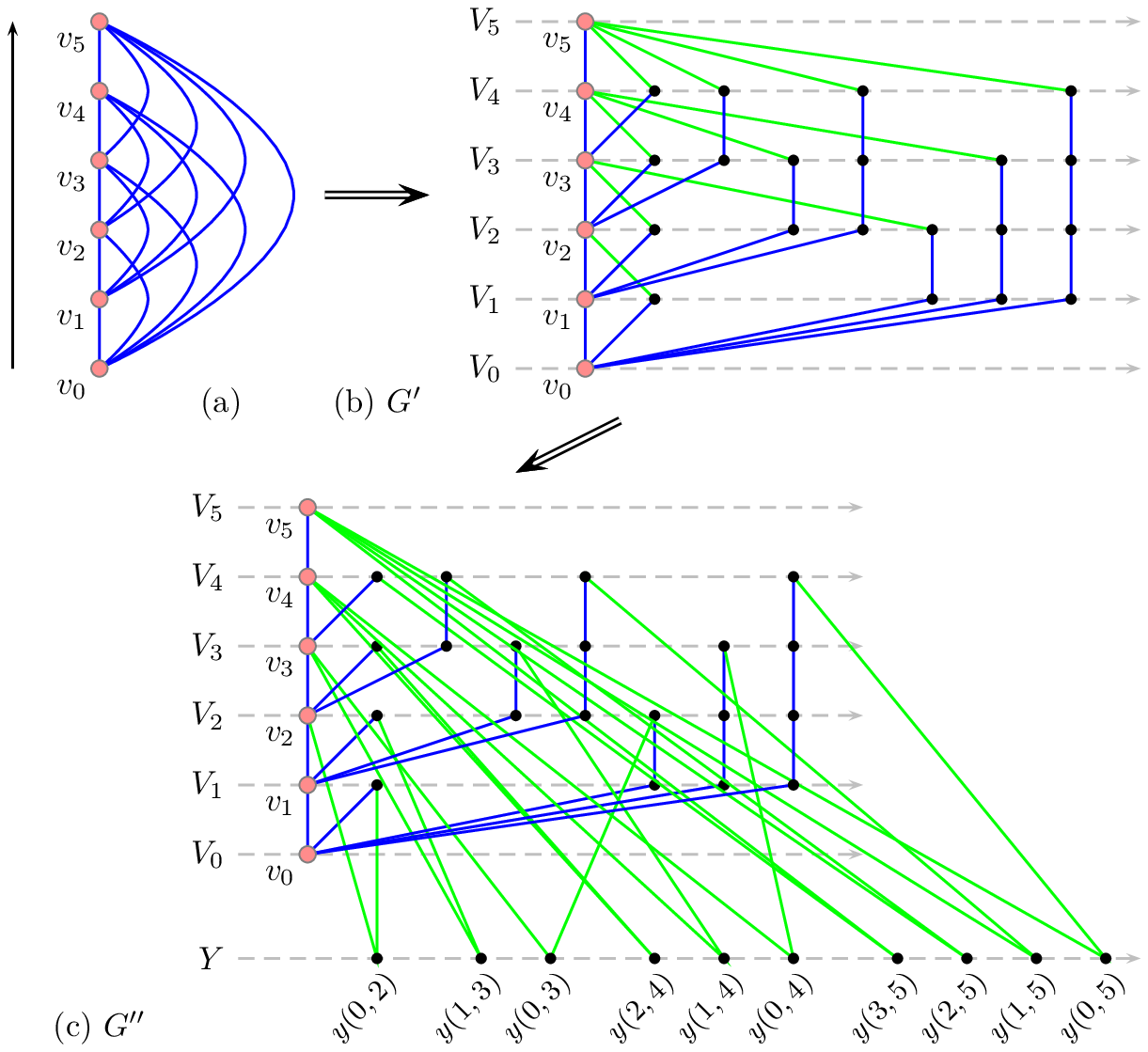}}{Track layouts in the proof of \thmref{FourTrackSub}.}

Order the green arcs $g(i,j)$ by non-decreasing $j$, breaking ties by decreasing $i$. Let $\pi$ be the total order obtained. Consider two green arcs $g_1$ and $g_2$ with $g_1<_\pi g_2$. Observe that if $g_1$ has an endpoint $p$ in the same track as an endpoint $q$ of $g_2$, then $p\leq q$ in this track. Call this property ($\star$). 

For every arc $\arc{v_iv_j}\in A(G)$ with $j\geq i+2$, subdivide the green arc $g(i,j)$ once to obtain a subdivision $G''$ of $G$ with at most $b$ division vertices per arc. Denote the division vertex by $y(i,j)$. As illustrated in \figref{FourTrack}(c), place all the vertices $y(i,j)$ in a new track $Y$ ordered by $\pi$. By property ($\star$), we obtain an $(n+1)$-track layout of $G''$ with no X-crossing. 

Now $\{V_1,V_2,\dots,V_n\}$ is an upward $n$-track layout of the subgraph of $G''$ consisting of the blue arcs, all of which have span one. By \lemref{Wrap}(c), $\{V_1,V_2,\dots,V_n\}$ can be wrapped into an upward $3$-track layout $\{X_0,X_1,X_2\}$, where $X_i=(V_i,V_{i+3},V_{i+4},\dots)$. Observe that property ($\star$) is maintained. Hence $\{X_0,X_1,X_2,Y\}$ is a $4$-track layout of $G''$. To see that $\{X_0,X_1,X_2,Y\}$ is upward, think of the $y(i,j)$ vertices as being in the middle of the corresponding green arcs in \figref{FourTrack}(b); then all the arcs in $(G'')^+$ point up or to the right. Thus $(G'')^+$ is acyclic, and $\{X_0,X_1,X_2,Y\}$ is an upward $4$-track layout of $G''$.
\end{proof}

%Note that the green arcs induce a forest of stars, each rooted at some vertex $v_j$ where $3\leq j\leq n$. We proceed to the second stage of the construction. Let $S_i$, $i\geq 3$, denote the green star rooted at $v_i$. In $T'$ each vertex of star $S_p$ comes before any vertex of star $S_q$, for all $3\leq p<q\leq n$, so stars are `nicely' ordered, as illustrated in FIGURE 2. Subdivide each green arc of star $S_i$ once and denote the resulting division vertex by $x(p,i,i)$, $p\leq i-2$. The resulting subdivision $G''$ or $G$ has at most $b$ division vertices per edge. We obtain $4$-track layout $T''$ of $G''$ from the $(2,3)$-track layout $T'$ by assigning all new division vertices $x(p,i,i)$, $3\leq i\leq n$, $p\leq i-2$, to the new track, track four. Their order in the track is determined first by the order of stars, breaking ties by the order of the leafs in stars, as illustrated in \figref{FIGURE 3}. There are no X-crossings thus $T''$ is a $4$-track layout of $G''$. $T''$ is upward since .... ({\tt How do we explain why DAVID? I see it is upward geometrically buy adding(in my head) all the extra arcs between green division vertices in Figure 1 and then seeing geometrically that they are all point up}). 

%%%%%%%%%%%%%%%%%%%%%%%%%%%%%%%%%%%%%%%%%%%

\twothmref{TwoQueueSub}{FourTrackSub} are best possible in the sense that every subdivision $H$ of a non-planar dag has $\uqn{H}\geq2$ and $\utn{H}\geq4$ (since $1$-queue graphs and $3$-track graphs are planar, and subdividing edges preserves planarity). In the undirected case, \citet{DujWoo-DMTCS05} proved that planarity characterises those graphs with $1$-queue or $3$-track subdivisions. That is, a graph $G$ is planar if and only if $G$ has a $1$-queue subdivision if and only if $G$ has a $3$-track subdivision. While we have not found such a characterisation for upward layouts, the proof of sufficiency by \citet[Lemma~32]{DujWoo-DMTCS05} generalises as follows. A dag $G$ is \emph{upward planar} if $G$ has a crossing-free drawing in the plane, such that every arc $\arc{vw}\in A(G)$ is represented by a \y-monotone curve with $\y(v)<\y(w)$. 

\begin{lemma}
Every upward planar dag $G$ has a subdivision that admits an upward $1$-queue layout and an upward $3$-track layout.
\end{lemma}

\begin{proof}
Given an upward planar drawing of $G$, draw a horizontal line through every  vertex. Now subdivide every edge whenever it crosses such a horizontal line. We obtain an upward track layout with span one, which by \lemref{UpwardOneQueue} has an upward $1$-queue layout and an upward $3$-track layout.
\end{proof}

%%%%%%%%%%%%%%%%%%%%%%%%%%%%%%%%%%%%%%%%%%%%%%%%%%%%%%%%%%%%%%%%%%%%%%%%%%%%%%
\subsection{Upward Polyline Drawings}
%%%%%%%%%%%%%%%%%%%%%%%%%%%%%%%%%%%%%%%%%%%%%%%%%%%%%%%%%%%%%%%%%%%%%%%%%%%%%%

A \emph{$b$-bend 3D drawing} of a graph $G$ is a \emph{3D drawing} of a subdivision of $G$ with at most $b$ division vertices per edge. \lemref{FourTrackDrawing} and \thmref{FourTrackSub} imply the following result.

\begin{corollary}
Let $G$ be an $n$-vertex $m$-arc dag with directed bandwidth $b$. Then $G$ has an upward $b$-bend $2\times 2\times(n+bm)$ drawing.\qed
\end{corollary}

The following theorem is a generalisation of a result by \citet{DujWoo-DMTCS05}.

\begin{theorem}
\thmlabel{TwoBendDrawing}
Every $n$-vertex dag $G$ with upward queue-number $\uqn{G}\leq k$ has an upward 2-bend $2k\times 2\times 2n$ drawing.
\end{theorem}

\begin{proof}
Let $(v_1,v_2,\dots,v_n)$ be a topological vertex ordering of $G$ that admits an upward $k$-queue layout. Number the queues $0,1,\dots,k-1$. Put each vertex $v_i$ at $(0,0,2i)$. Draw each arc \arc{v_iv_{i+1}} straight. For all $j\geq i+2$, draw each arc \arc{v_iv_j} in queue $\ell$ with the $2$-bend polyline \begin{equation*}
(0,0,2i)\rightarrow(2\ell,1,i+j)\rightarrow(2\ell+1,1,i+j+1)\rightarrow(0,0,2j).
\end{equation*} 
Since $2i<i+j<i+j+1<2j$, the drawing is upward. \citet[Theorem~26]{DujWoo-DMTCS05} proved that no two arcs cross.
%Two arcs in distinct queues do not cross, since their projections into the $XY$-plane only intersect at $(0,0)$. Consider two arcs $v_iv_j$ and $v_av_b$ in queue $\ell$. Without loss of generality $i\leq a$, $j\leq b$, and $i+j<a+b$. By again considering the projection into the $XY$-plane, a crossing between $v_iv_j$ and $v_av_b$ must involve the first segment of both edges, the second segment of both edges, or the third segment of both edges. The first bends are at $(2\ell,1,i+j)$ and $(2\ell,1,a+b)$. Since $i\leq a$ and $i+j<a+b$, the first segments do not cross. The second bends are at $(2\ell+1,1,i+j+1)$ and $(2\ell+1,1,a+b+1)$. Since $i+j<a+b$ and $i+j+1<a+b+1$, the second segments do not cross. 
\end{proof}

It follows from results of \citet{HR-SJC92} that the upward queue-number of every $n$-vertex graph is at most $\floor{\frac{n}{2}}$ (which is tight for the complete dag). Thus we have the following corollary of \thmref{TwoBendDrawing}.

\begin{corollary}
\corlabel{TwoBendDrawing}
Every $n$-vertex dag $G$ has an upward $2$-bend $n\times 2\times 2n$ drawing with volume $4n^2$.\qed
\end{corollary}

%%%%%%%%%%%%%%%%%%%%%%%%%%%%%%%%%%%%%%%%%%%%%%%%%%%%%%%%%%%%%%%%%%%%%%%%%%%%%%%%
%\bibliographystyle{myNatbibStyle}
%\bibliography{myBibliography,myConferences}
%%%%%%%%%%%%%%%%%%%%%%%%%%%%%%%%%%%%%%%%%%%%%%%%%%%%%%%%%%%%%%%%%%%%%%%%%%%%%%%%

\def\soft#1{\leavevmode\setbox0=\hbox{h}\dimen7=\ht0\advance \dimen7
  by-1ex\relax\if t#1\relax\rlap{\raise.6\dimen7
  \hbox{\kern.3ex\char'47}}#1\relax\else\if T#1\relax
  \rlap{\raise.5\dimen7\hbox{\kern1.3ex\char'47}}#1\relax \else\if
  d#1\relax\rlap{\raise.5\dimen7\hbox{\kern.9ex \char'47}}#1\relax\else\if
  D#1\relax\rlap{\raise.5\dimen7 \hbox{\kern1.4ex\char'47}}#1\relax\else\if
  l#1\relax \rlap{\raise.5\dimen7\hbox{\kern.4ex\char'47}}#1\relax \else\if
  L#1\relax\rlap{\raise.5\dimen7\hbox{\kern.7ex
  \char'47}}#1\relax\else\message{accent \string\soft \space #1 not
  defined!}#1\relax\fi\fi\fi\fi\fi\fi} \def\cprime{$'$}

\end{document}